\documentclass[12pt]{article}
\usepackage[utf8]{inputenc}
\usepackage{amssymb,amsmath,amsfonts,amsthm,amscd,latexsym,indentfirst,verbatim}
\usepackage{geometry}
\usepackage{geometry}
\def\Imm{\textrm{Im}}
\def\Aut{\textrm{Aut}}
\def\Span{\textrm{Span}}
\def\Spec{\textrm{Spec}}
\def\id{\textrm{id}}
\def\gcd{\textrm{gcd}}

\begin{document}

\begin{center}
{\Large Monomial Rota---Baxter operators on free commutative non-unital algebra}

Vsevolod Gubarev
\end{center}



\begin{abstract}
A Rota---Baxter operator defined on the polynomial algebra
is called monomial if it maps each monomial to a monomial with some coefficient.
We classify monomial Rota---Baxter operators defined on
the algebra of polynomials in one variable without constant term.
We also describe injective monomial Rota---Baxter operators of nonzero weight
on the algebra of polynomials in several variables without constant term.

\end{abstract}

\section{Introduction}

A linear operator $R$ defined on a (non-asssociative, in general) algebra~$A$
over a field~$\Bbbk$ is called a Rota---Baxter operator (RB-operator), if the following relation
\begin{equation}\label{RB}
R(x)R(y) = R(R(x)y + xR(y) + \lambda xy)
\end{equation}
holds for all $x,y\in A$.
Here $\lambda$ is a fixed constant from $\Bbbk$ called a weight of $R$.
An algebra equipped with a Rota---Baxter operator is called a Rota---Baxter algebra.

G. Baxter introduced the notion of a Rota---Baxter operator in 1960~\cite{Baxter}
as a natural generalization of integration by parts formula for the integral operator.
Further, many authors including G.-C. Rota, P. Cartier, L. Guo have studied RB-operators,
see details in~\cite{Unital,GuoMonograph}.
There are deep connections of Rota---Baxter algebras with mathematical physics, number theory,
operad theory, combinatorics, and others.

One of the directions in the theory of Rota---Baxter operators
is to study  RB-ope\-rators on polynomials~\cite{Monom2,Monom}
and power series~\cite{GuoLiu,Mazurek}.
Since both algebras are infi\-nite-dimensional,
the classification of all RB-operators on them seems to be a hard problem.
Thus, RB-operators of special kind were introduced into consideration.
One of such special RB-operators is the class of monomial RB-operators~\cite{Monom2},
i.e., such RB-operators that map every monomial to some monomial with (maybe, zero) coefficient.
L. Guo, M.~Rosenkranz, and S.H. Zheng in 2015 described all injective
monomial RB-operators of weight zero on $\Bbbk[x]$.
In 2016~\cite{Monom}, H. Yu classified all monomial RB-operators of any weight on $\Bbbk[x]$.

In the study of RB-operators, the presence of the multiplicative identity (unit)
in the algebra is a crucial condition. As shown in \cite{Unital}, RB-operators
on a unital algebra subject to various additional restrictions in comparison with
non-unital algebras.

In the present paper, we classify all monomial RB-operators of weight zero (Theorem~1)
and nonzero weight (Theorem~2) on $\Bbbk_0[x]$, the free commutative non-unital
algebra generated by $x$. We obtain the complete classification of H.~Yu as
a corollary (Corollaries~2 and 3).
Further, we describe injective monomial RB-operators of nonzero weight on both
$\Bbbk_0[X]$ (Theorem~3) and $\Bbbk[X]$ (Corollary~4).
Connected to the last results, a partial grading on a finite-dimensional algebra
by the spectrum of its RB-operator (Propositions~2 and~3) is stated.
The analogous grading in the case of derivations and automorphisms is well-known~\cite{Jacobson}.
As examples of such partial gradings we consider monomial RB-operators
on the quotient of $\Bbbk_0[x]$ by the ideal generated by $x^{N+1}$ for some $N>0$
(Examples~7 and~9).

Excepting the last section, we assume that a ground field $\Bbbk$ has characteristic zero.

\section{Preliminaries}

Let us start with some basic properties of Rota---Baxter operators.

Trivial RB-operators of weight $\lambda$ are zero operator and $-\lambda\id$.

{\bf Lemma 1}~\cite{BGP,GuoMonograph}.
Let $A$ be an algebra and let $P$ be an RB-operator of weight $\lambda$ on $A$.

a) The operator $\lambda^{-1}P$ is an RB-operator of weight 1 provided that $\lambda\neq0$,

b) Given an automorphism $\psi\in\Aut(A)$, the operator
$\psi^{-1}P\psi$ is an RB-operator of weight~$\lambda$ on $A$.

{\bf Lemma 2}~\cite{GuoMonograph}.
Let an algebra $A$ split as a vector space into a direct sum of two subalgebras
$A_1$ and $A_2$. An operator $P$ defined as $P(a_1+a_2) = -\lambda a_2$,
$a_1\in A_1$, $a_2\in A_2$, is an RB-operator of weight~$\lambda$ on $A$.

We call an RB-operator from Lemma~1 as {\it splitting RB-operator} with subalgebras $A_1$ and $A_2$.

{\bf Lemma 3}~\cite{BGP,GuoMonograph}.
Let $A$ be a unital algebra and let $P$ be an RB-operator of weight $\lambda$ on $A$.

a) If $\lambda\neq0$ and $P(1)\in \Bbbk$, then $P(1)\in\{0,-\lambda\}$ and $P$ is splitting.

b) If $\lambda = 0$ and $P(x) \in \Bbbk$, then $P(x) = 0$.

Let $R$ be an RB-operator of weight $\lambda$ on an algebra~$A$.
Then $\Imm(R)$ is a subalgebra of~$A$.
If $\lambda\neq0$, then $\ker(R)$ is also a subalgebra of $A$.
If $\lambda = 0$, then $\ker(R)$ is an $\Imm(R)$-module.

For $\lambda = 0$, the relation
\begin{multline}\label{LongRB}
R(x_1)R(x_2)\ldots R(x_k)
 = R(x_1R(x_2)\ldots R(x_k)+R(x_1)x_2R(x_3)\ldots R(x_k) \\
 +\ldots+R(x_1)R(x_2)\ldots R(x_{k-1})x_k)
\end{multline}
holds in an associative RB-algebra as a direct consequence of~\eqref{RB}.
In particular, when RB-algebra is commutative, $\lambda = 0$, and
$x = x_1 = x_2 = \ldots = x_k$,
we have
\begin{equation}\label{LongRBEqual}
(R(x))^k = k R(x(R(x))^{k-1}).
\end{equation}

From now on, we focus on RB-operators defined on either polynomial algebra or
on free commutative non-unital algebras, i.e., polynomials without constant term.

{\bf Example 1}~\cite{GuoMonograph}.
The linear map $P$ defined on $\Bbbk[x]$ by the formula
$J_a(x^n) = \frac{x^{n+1}-a^{n+1}}{n+1}$, where $a$ is a fixed element from $\Bbbk$,
is an RB-operator of weight zero on $\Bbbk[x]$.

Given an algebra $A$ and an element $r\in A$,
denote by $l_r$ the linear operator on~$A$ defined by $l_r(x) = rx$, $x\in A$.

{\bf Example 2}~\cite{Monom2}.
Let $A$ be a commutative algebra, $r\in A$, and let $R$ be an RB-operator of weight zero on~$A$.
Then the linear map $R\circ l_r$ is again an RB-operator
of weight zero on~$A$.
Here $\circ$ denotes the composition of operators.

Thus, the linear operator $P = J_a\circ l_{x^k}$ acting on $\Bbbk[x]$ as
$P(x^n) = \frac{x^{n+k+1}-a^{n+k+1}}{n+k+1}$ is an RB-operator of weight zero.

Given an algebra $A$, it is known that every solution
$r = \sum a_i\otimes b_i\in A\otimes A$
of the associative Yang---Baxter equation of weight~$\lambda$ (AYBE)~\cite{FardThesis,Zhelyabin}
\begin{equation}\label{wAYBE}
r_{13}r_{12} - r_{12}r_{23} + r_{23}r_{13} = \lambda r_{13}
\end{equation}
gives rise to an RB-operator of weight~$-\lambda$ on $A$~\cite{Aguiar,FardThesis} defined by the formula
\begin{equation}\label{AguiarsFormula}
R(x) = \sum\limits a_i x b_i.
\end{equation}
In~\eqref{wAYBE}, $r_{12} = \sum\limits a_i\otimes b_i\otimes 1$,
$r_{13} = \sum\limits a_i\otimes 1\otimes b_i$, and
$r_{23} = \sum\limits 1\otimes a_i\otimes b_i$.

{\bf Proposition 1}.
The only nonzero solution of AYBE of weight~$\lambda\neq0$ on $\Bbbk[x]$
is $\lambda(1\otimes 1)$.

{\sc Proof}.
Let $r = \sum\limits_{i,j\geq 0}\alpha_{ij}x^i\otimes x^j$
be a solution of~\eqref{wAYBE}, so we have
\begin{equation}\label{AYBEk[x]}
\sum\limits_{i,j,k,l\geq0}\alpha_{ij}\alpha_{kl}
(x^{i+k}\otimes x^l\otimes x^j - x^i\otimes x^{j+k}\otimes x^l
+ x^i\otimes x^k\otimes x^{j+l})
- \lambda\sum\limits_{i,j\geq0} x^i\otimes 1\otimes x^j = 0.
\end{equation}
Consider maximal $N$ such that $\alpha_{Nj}\neq0$ for some $j$.
If $N>0$, then the left-hand side of~\eqref{AYBEk[x]} is nonzero
because of the summand
$\alpha_{Nj}^2 x^{2N}\otimes x^j\otimes x^j$
from the first sum.
Analogously we may consider maximal $M$
such that $\alpha_{0M}\neq0$.
So, the only possible solution is a tensor $q(1\otimes 1)$.
It is easy to get that either $q = \lambda$ or $q = 0$.

{\bf Corollary 1}.
The only nonzero solution of AYBE of weight~$\lambda\neq 0$ on $\Bbbk[X]$ is $\lambda(1\otimes 1)$.

Given a nonempty set $X$, by $\Bbbk_0[X]$ we denote the free commutative algebra generated by $X$.
By the formula~\eqref{AguiarsFormula}, we get only trivial RB-operators on $\Bbbk[X]$ and $\Bbbk_0[X]$.

\section{Monomial RB-operators of weight zero on $\Bbbk_0[x]$}

A linear operator $R$ defined on $\Bbbk[x]$ ($\Bbbk_0[x]$) is called monomial if for all $n$
we have $R(x^n) = \alpha_n x^{t_n}$ for some $\alpha_n\in \Bbbk$ and $t_n\in \mathbb{N}$
($t_n\in \mathbb{N}_{>0}$).

{\bf Theorem 1}.
Given a nonzero monomial RB-operator $R$ on $\Bbbk_0[x]$,
there exist positive $m\in\mathbb{N}$, nonnegative
$p_0,\ldots,p_{m-1}\in\mathbb{N}$ and some
$q_0,\ldots,q_{m-1}\in \Bbbk$ such that
$p_i = 0$ if and only if $q_i = 0$,
and $R$ is defined by the following formula
\begin{equation}\label{Monomialk[x]ZeroNonunital}
R(x^{ma+b})
 = q_b\frac{x^{m(a+p_b)}}{m(a+p_b)},
\end{equation}
where $a\in\mathbb{N}$ and $0<b\leq m$.

{\sc Proof}.
Since $R$ is nonzero, at least one monomial lies in $\Imm(R)$.
So, $\Imm(R)$ is an infinite-dimensional subalgebra of $\Bbbk[x]$.
Define $m = \gcd(t\mid x^t\in\Imm(R))$.
We have
$\Imm(R) \supset \Span\{x^{mk}\mid k\geq N\}$ for some natural $N$.

{\bf Lemma 4}.
Let $0<b\leq m$. Suppose that $x^{ma + b}\in\ker(R)$ for some $a\geq0$.
Then $x^{mc + b}\in\ker(R)$ for all $c\geq0$.

{\sc Proof of Lemma~4}.
Since by~\eqref{RB} $\ker(R)$ is an $\Imm(R)$-module,
$x^{mc + b}\in\ker(R)$ for all $c\geq a + N$.
Assume there exists $c$ such that
$R(x^{mc+b}) = \alpha x^{mt}$ and $\alpha\neq0$.
Then by~\eqref{LongRBEqual}, we get
$$
\alpha^k x^{mtk}
 = (R(x^{mc+b}))^{k}
 = k R(x^{mc+b}(R(x^{mc+b}))^{k-1})
 = k\alpha^{k-1}R(x^{m(c+t(k-1))+b}).
$$
We obtain a contradiction when $k\geq a + N + 1$.
Lemma~4 is proved.

{\bf Lemma 5}.
Let $0<b\leq m$. Suppose that $R(x^b) = \alpha_0 x^{mp_b}$ with $\alpha_0\neq0$.
Then $R(x^{ma + b}) = \alpha_a x^{m(a+p_b)}$ for all $a\geq0$.
Here $\alpha_a$ are some nonzero elements from $\Bbbk$.

{\sc Proof of Lemma~5}.
First, we prove Lemma~5 for all $a\geq N$.
Since $x^{ma}\in\Imm(R)$, we may find $k\geq0$
such that $R(x^k) = \gamma x^{ma}$ for some $\gamma\neq0$.
Then
$$
\alpha_0\gamma x^{m(a+p_b)}
 = R(x^b)R(x^k)
 = \alpha_0 R(x^{mp_b + k}) + \gamma R(x^{ma+b}).
$$
By Lemma~4, $R(x^{ma+b})\neq0$.
Since $R$ is monomial, we have
$R(x^{ma+b}) = \alpha_a x^{m(a+p_b)}$ for some $\alpha_a\in \Bbbk$.

Now, consider $0<a<N$.
Suppose $R(x^{ma+b}) = \beta x^{mt}$ for some $t>0$.
By Lemma~4, $\beta\neq0$. Then
$$
\alpha_0\beta x^{m(t+p_b)}
 = R(x^b)R(x^{ma+b})
 = \alpha_0 R(x^{m(a+p_b)+b}) + \beta R(x^{mt+b}).
$$
If $t = a+p_b$, we are done.
If $t\neq a+p_b$, then the monomials
$x^{m(a+p_b)+b}$ and $x^{mt+b}$ have proportional images under $R$,
it means there exist $k<l$ and $\delta\neq0$ such that
$x^{mk+b} + \delta x^{ml+b}\in\ker(R)$.
Multiplying, if necessary, this element by $x^{mc}\in\Imm(R)$ for $c\geq N$,
we may assume  $k\geq N$.
But for such degrees we have already proved that
$R(x^{mk+b} + \delta x^{ml+b})
 = \alpha_k x^{m(k+p_b)} + \alpha_l\delta x^{m(l+p_b)}\neq0$,
a~contradiction. Lemma~5 is proved.

We may rewrite each positive natural number~$n$ as $n = ma + b$ for $m\geq0$
and $0<b\leq m$. Define $\alpha_{a,b}\in \Bbbk$ so that
$R(x^{ma+b}) = \alpha_{a,b}x^{m(a+p_b)}$.
If $x^b\in\ker(R)$, then $\alpha_{0,b} = p_b = 0$.
Otherwise, $\alpha_{0,b},p_b\neq0$.
Thus, the identity~\eqref{RB}
\begin{multline*}
\alpha_{a,b}\alpha_{c,b'} x^{m(a+c+p_b+p_{b'})}
 = R(x^{ma+b})R(x^{mc+b'}) \\
 = \alpha_{a,b}R(x^{m(a+c+p_b)+b'})
  + \alpha_{c,b'}R(x^{m(a+c+p_{b'})+b}) \\
 = (\alpha_{a,b}\alpha_{a+c+p_b,p_{b'}}
 +\alpha_{c,b'}\alpha_{a+c+p_{b'},p_b})x^{m(a+c+p_b+p_{b'})}
 \end{multline*}
is equivalent to the following one:
\begin{equation}\label{MonomialZeroCoefficients}
\alpha_{a,b}\alpha_{c,b'}
 = \alpha_{a,b}\alpha_{a+c+p_b,p_{b'}}
 +\alpha_{c,b'}\alpha_{a+c+p_{b'},p_b}.
\end{equation}
By~Lemma~4, it is reasonable to study only the cases when $x^b,x^{b'}\not\in\ker(R)$.
Put $b = b'$ and $a = c$ in~\eqref{MonomialZeroCoefficients},
then $\alpha_{a,b} = 2\alpha_{2a+pb,b}$.
This relation can be rewritten as
$\gamma_{a,b} = \gamma_{2a+p_b,b}$, where
$\gamma_{a,b} = (a+p_b)\alpha_{a,b}$.
For $a = 0$, we get $\gamma_{0,b} = \gamma_{p_b,b}$.
Let us prove that $\gamma_{np_b,b} = \gamma_{0,b}$ by induction on $n$.
We have got the base case $n=1$.
Suppose this equality holds for all natural numbers
less or equal to~$n$. The relation~\eqref{MonomialZeroCoefficients},
for $a = 0$, $c = np_b$, and $b' = b$, turns into
$$
\alpha_{0,b}\alpha_{np_b,b}
 = \alpha_{(n+1)p_b,b}(\alpha_{0,b}+\alpha_{np_b,b}).
$$
Applying the induction hypothesis, we obtain
\begin{multline*}
\gamma_{(n+1)p_b,b}
 = (n+2)p_b\alpha_{(n+1)p_b,b}
 = (n+2)p_b\frac{\alpha_{0,b}\alpha_{np_b,b}}{\alpha_{0,b}+\alpha_{np_b,b}} \\
 = (n+2)p_b \frac{\gamma_{0,b}\gamma_{np_b,b}}{(n+1)p_b^2}\cdot
 \frac{1}{\gamma_{0,b}(1/p_b+1/((n+1)p_b)}
 = \gamma_{np_b,b}
 = \gamma_{0,b}.
\end{multline*}

Further, let us show that $\gamma_{a,b} = \gamma_{0,b}$ for all $a$.
Indeed, the relation~\eqref{LongRBEqual} applied for $x^{ma+b}$ leads to
$$
\alpha_{a,b}^k
 = k\alpha_{a,b}^{k-1}\alpha_{ka+(k-1)p_b,b}
$$
or
$\gamma_{a,b} =\gamma_{ka+(k-1)p_b,b}$.
Choosing $k = p_b$, we conclude by the above property
that $\gamma_{a,b} = \gamma_{0,b}$.

Finally, define $q_b = m\gamma_{0,b}$ for all $0<b\leq m$.
We see that $q_b = 0$ if and only if $p_b = 0$ by the definition.
It is easy to check that the linear operator $R$ obtained
is an RB-operator of weight~0 on $\Bbbk_0[x]$.
Theorem is proved.

{\bf Example 3}.
The linear map $R$ on $\Bbbk_0[x]$ defined as
$R(x^n) = \frac{x^n}{n}$ is an RB-operator of weight zero.
It is easy to see that $R$ is invertible and so $d = R^{-1}$
is a derivation on $\Bbbk_0[x]$ with $d(x^n) = nx^n$.
Hence, the restriction of the RB-operator $R$
from Example~1 to $\Bbbk_0[x]$ equals $d^{-1}\circ x$.

{\bf Corollary 2}~\cite{Monom}.
Given a nonzero monomial RB-operator $R$ on $\Bbbk[x]$,
there exist positive $m\in\mathbb{N}$,
nonnegative $p_1,\ldots,p_m\in\mathbb{N}$
and some $q_1,\ldots,q_m\in \Bbbk$ such that
$p_i = 0$ if and only if $q_i = 0$, and $R$ is defined by
$$
R(x^{ma+b})
 = q_b\frac{x^{m(a+p_b)}}{m(a+p_b)},
$$
where $a\in\mathbb{N}$ and $0\leq b<m$.

{\sc Proof}.
By Lemma~3b, $1\not\in \Imm(R)$, so $\langle x\rangle$ is $R$-invariant.
We may apply Theorem~1. It remains to extend $R$ from $\langle x\rangle$ to~$\Bbbk[x]$.
It is enough to define $R(x^0) = \alpha_{0,0}x^{mp_0}$
instead of the previous $R(x^m) = \alpha_{0,m}x^{mp_m}$
and we analogously prove that $\gamma_{a,m} = \gamma_{0,m}$ for all $a\geq 0$.

Let us give a few remarks about the classification of monomial RB-operators of weight zero on $\Bbbk[x]$.

{\bf Remark 1}.
When $x^m\in\ker(R)$, the RB-operator $R$ is defined by Lemma 5.1b from~\cite{Unital}
with the decomposition $\Bbbk[x] = B\oplus C$ (of vector spaces) for
$B = \Span\{x^{am+b}\mid p_b = 0\}$ and
$C = \Span\{x^{am+b}\mid p_b\neq0\}$.

{\bf Remark 2}.
If $x^m\not\in\ker(R)$, then we may consider the restriction of $R$
to the subalgebra $\Span\{x^{mk}\mid k\geq0\}$.
Then, up to a scalar multiple, $R$ is defined on $\Bbbk[y]$
for $y = x^m$ by the formula
$R(y^a) = \frac{y^{a+p_0}}{a+p_0}$.
In other words, $R = J_0\circ l_{y^{p_0-1}}$.
So, $R$ is analytically modeled RB-operator
on $\Bbbk[x]$ in the terminology of~\cite{Monom2}.

\section{Monomial RB-operators of nonzero weight on $\Bbbk_0[x]$}

Let $R$ be an RB-operator of nonzero weight $\lambda$ on an algebra~$A$.
By Lemma~1a, we may assume that $\lambda = 1$.

{\bf Theorem 2}.
Let $R$ be a monomial RB-operator of weight~1 on $\Bbbk_0[x]$.
Then there exists $\alpha\in \Bbbk$ such that $(\alpha+1)^n \neq \alpha^n$ for all $n\geq1$ and
\begin{equation}\label{RBNonzeroNonunitary}
R(x^n) = \frac{\alpha^n}{(\alpha+1)^n-\alpha^n}x^n.
\end{equation}

{\sc Proof}.
Trivial RB-operators on $\Bbbk_0[x]$ are monomial ones
and they correspond to the cases $\alpha=0$ and $\alpha=-1$ in~\eqref{RBNonzeroNonunitary}.
Suppose that $R$ is a nontrivial monomial RB-operator of weight~1 on $\Bbbk_0[x]$.

\underline{{\sc Case 1}: $\ker(R) = (0)$}.
Suppose that there exists a monomial $x^n$ such that $R(x^n) = \alpha x^m$ for $n\neq m$.
Since $R$ is injective, $\alpha\neq0$. Then
$$
\alpha^2 x^{2m}
 = R(x^n)R(x^n)
 = 2\alpha R(x^{n+m}) + R(x^{2n}).
$$
Since $n+m\neq 2n$, the kernel of $R$ is nonzero, a contradiction.
Thus, $R(x^n) = \alpha_n x^n$, $\alpha_n\in \Bbbk$, for all $n\geq1$.

We have $R(x) = \alpha_1 x$, denote $\alpha = \alpha_1$.
Let us prove the formula~\eqref{RBNonzeroNonunitary} by induction on~$n$.
For $n=1$, we are done. Suppose that~\eqref{RBNonzeroNonunitary}
is proved for all numbers less than $n$. Then
\begin{gather*}
R(x)R(x^{n-1})
 = \alpha x\cdot\frac{\alpha^{n-1}}{(\alpha+1)^{n-1} -\alpha^{n-1} }x^{n-1}
 = \frac{\alpha^n}{(\alpha+1)^{n-1}-\alpha^{n-1}}x^n, \\
R(x)R(x^{n-1})
 = R(R(x)x^{n-1} + xR(x^{n-1}) + x^n)
 = \left(\alpha + \frac{\alpha^{n-1}}{(\alpha+1)^{n-1}-\alpha^{n-1}} + 1\right)R(x^n).
\end{gather*}
Thus,
$$
R(x^n)
 = \frac{\alpha^n}{(\alpha+1)^{n-1}-\alpha^{n-1}}
 \cdot \frac{(\alpha+1)^{n-1}-\alpha^{n-1}}{\alpha^{n-1}
 + (\alpha+1)^n -\alpha^n - \alpha^{n-1} }x^n
 = \frac{\alpha^n}{(\alpha+1)^n-\alpha^n}x^n,
$$
as required.

\underline{{\sc Case 2}: $\ker(R)\neq(0)$}.
Since $R$ is nontrivial, $\Imm(R)\neq(0)$ as well.
Define $m = \gcd(t\mid x^t\in\Imm(R))$. It is easy to get that
$$
\Imm(R) \supset \Span\{x^{mk}\mid k\geq N\}
$$
for some natural $N$.

\underline{{\sc Case 2a}: $\ker(R)$ contains a monomial $x^k$}.
Then $x^{mkN}\in \ker(R)\cap \Imm(R)$, since
$\ker(R)$ and $\Imm(R)$ are subalgebras of $\Bbbk_0[x]$.
Denote $t = kN$ and consider $x^s$ such that $R(x^s) = ax^{mt}$, $a\neq0$.
It is known that $\ker(R)$ is an ideal in $\Imm(R+\id)$, i.e.,
$x^{ps+qmt}\in\ker(R)$ for all $p,q\geq1$. Then
$R(x^s)R(x^s) = a^2 x^{2mt} = R(x^{2s})$
and, analogously, $R(x^{smt}) = a^{mt}x^{m^2t^2}$.
We have $x^{smt}\in\ker(R)$, so $a = 0$, a contradiction.

\underline{{\sc Case 2b}: $\ker(R)$ does not contain a monomial}.
Since $R$ is monomial, $\ker(R)$ is non\-zero only if there exist
$p\neq s$ such that $R(x^p) = ax^{mt}$ and $R(x^s) = bx^{mt}$
for some $t\geq1$ and $a,b\neq0$.
We have
$$
a^2 x^{2mt}
 = R(x^p)R(x^p)
 = 2aR(x^{p+mt}) + R(x^{2p}),
$$
so $R(x^{2p}) = \alpha_{2p}x^{2mt}$.
By induction, we may prove that
$R(x^{pk}) = \alpha_{pk}x^{mtk}$ for all $k\geq1$.
Analogously, $R(x^{sk}) = \alpha_{sk}x^{mtk}$ for all $k\geq1$.
Let us note that simultaneously
$R(x^{ps}) = \alpha_{ps}x^{mts}$
and $R(x^{ps}) = \alpha_{ps}x^{mtp}$.
We arrive at a contradiction since $p\neq s$.
Theorem is proved.

{\bf Corollary 3}~\cite{Monom}.
Up to conjugation with an automorphism of $\Bbbk[x]$, each nontrivial
monomial RB-operator on $\Bbbk[x]$ of nonzero weight is splitting with
subalgebras $\Bbbk$ and $\langle x\rangle$.

{\sc Proof}.
Let $R$ be a nontrivial monomial RB-operator of weight~1 on $\Bbbk[x]$.
Let $R(1) = \alpha x^k$ and suppose that $k>0$ and $\alpha\neq0$.
Then by~\eqref{RB},
\begin{equation}\label{R(1)R(1)Nonzero}
2\alpha R(x^k)
 = 2 R^2(1)
 = R(1)R(1) - R(1)
 = \alpha^2 x^{2k} - \alpha x^k,
\end{equation}
a contradiction to the monomiality condition.

Thus, $R(1)\in \Bbbk$, i.e., $R(1) \in \{0,-1\}$, and $R$ is splitting by Lemma~3a.
Since $R$ is nontrivial, both $\ker(R)$ and $\Imm(R)$ are nonzero.
Hence, $\Imm(R)$ has a basis of monomials and $R$ acts on $\Imm(R)$ as the operator $-\id$.

{\sc Case 1}.
Suppose $\langle x\rangle = \Bbbk_0[x]$ is $R$-invariant.
By Theorem~2, $R(x^n) {=} \dfrac{\alpha^n}{(\alpha+1)^n-\alpha^n}x^n$, $n>0$,
for some $\alpha$ such that the denominator is nonzero for all $n$.
If $\alpha = 0$, then we have either trivial RB-operator (when $R(1) = 0$)
or the splitting one with subalgebras $\Bbbk$ and $\langle x\rangle$ (when $R(1) = -1$).

If $\alpha\neq0$, then since $R$ acts on $\Imm(R)$ as the operator $-\id$,
we get $\alpha = -1$ and $R(x^n) = -x^n$ for all $n>0$.
Again, we have either trivial RB-operator (when $R(1) = -1$)
or the splitting one with subalgebras $\Bbbk$ and $\langle x\rangle$ (when $R(1) = 0$).

{\sc Case 2}.
Suppose $\langle x\rangle = \Bbbk_0[x]$ is not $R$-invariant.
Let us show that $R(1) = -1$. Indeed, suppose that $R(1) = 0$
and there exists such $k>0$ that $R(x^k) = \alpha1$ for some $\alpha\in\Bbbk$.
Since $\ker(R)$ is an ideal in $\Imm(R+\id)$, we have $x^k\in \ker(R)$ and so, $\alpha = 0$.
Thus, $\langle x\rangle$ is $R$-invariant, a~contradiction.

Consider $R(x)$. If $R(x) = -x$, then $\ker(R+\id) = \Bbbk[x]$ and $R$ is trivial.
If $R(x) = 0$, then $R$ is splitting with subalgebras $\Bbbk$ and $\langle x\rangle$.
Assume that $R(x) = \alpha x^k$ for some $\alpha\neq0$ and $k>0$,
and moreover, $R(x)\neq -x$. From the last condition, we get $k>1$.
If $\ker(R)$ contains a monomial $x^t$, then $x^{tk}\in \ker(R)\cap \Imm(R) = (0)$,
a contradiction. Since $R(x^k) = - x^k$, we may repeat arguments of the proof of Theorem~2 to get a contradiction.

It remains to study the case when $R(x) = \alpha1$ for $\alpha\neq0$.
Let us assume for the convenience that the weight of the RB-operator $R$ equals~$-1$.
Then it is easy to show by induction on $n$ that
the equality $R(x) = \alpha1$ implies $R(x^n) = \alpha^n1$ for all $n\geq1$.
Define the automorphism $\varphi\colon x^n \to x^n/\alpha^n$ of $\Bbbk[x]$.
Then the RB-operator $R' = \varphi^{-1}R\varphi$ (we apply Lemma~1b)
acts on all monomials as $R'(x^n) = 1$
and $R'$ is splitting with subalgebras $\Imm(R') = \Bbbk$ and
$\ker(R') = \langle x-1\rangle$.
Define $\psi\in \Aut(\Bbbk[x])$ as $\psi\colon x\to x-1$.
Note that $\psi^{-1}R'\psi$ is the splitting RB-operator
with subalgebras $\Bbbk$ and $\langle x\rangle$.

\section{Monomial RB-operators of nonzero weight on $\Bbbk_0[X]$}

{\bf Example 4}.
Let $R$ be a splitting RB-operator on $\Bbbk_0[x,y]$
with subalgebras $\Bbbk_0[x]$ and $\langle y\rangle$.
Then $R$ is monomial.

{\bf Theorem 3}.\,Let $R$ be an injective monomial RB-operator of weight\,1 on $\Bbbk_0[x_1,\ldots\!,x_n]$,
then there exist nonzero
$\alpha_1,\ldots,\alpha_n\in \Bbbk$ such that
\begin{equation}\label{MonomialInjectiveOnk[X]}
R(x_1^{i_1}\ldots x_n^{i_n})
 = \frac{\alpha_1^{i_1}\ldots \alpha_n^{i_n}}
 {(\alpha_1+1)^{i_1}\ldots (\alpha_n+1)^{i_n}
 -\alpha_1^{i_1}\ldots \alpha_n^{i_n}}x_1^{i_1}\ldots x_n^{i_n}
\end{equation}
for all $i_1,\ldots,i_n\geq0$, $i_1^2+\ldots+i_n^2>0$.
Moreover, all denominators are nonzero.

{\sc Proof}.
First, suppose that there exists a monomial
$w = x_1^{i_1}\ldots x_n^{i_n}$ such that $R(w) = \alpha w'$ with $w'\neq w$.
Then
$$
\alpha^2 w'^2
 = R(w)R(w)
 = 2\alpha R(ww') + R(w^2).
$$
Since $ww'\neq w^2$, the kernel of $R$ is nonzero, a contradiction.

Define scalars $\alpha_i\in \Bbbk\setminus\{0\}$ such that
$R(x_i) = \alpha_i x_i$. Let us prove the formula~\eqref{MonomialInjectiveOnk[X]}
for a monomial~$w$ by induction on the degree $\deg(w)$.
Given a monomial $w$ from $\Bbbk_0[x_1,\ldots,x_n]$,
denote by $\alpha(w)$ the coefficient at $w$ in the right-hand side
of~\eqref{MonomialInjectiveOnk[X]}.
For $\deg(w) = 1$, \eqref{MonomialInjectiveOnk[X]}
follows from the definition of~$\alpha_i$.

Suppose that we have proved~\eqref{MonomialInjectiveOnk[X]}
for all monomials of degree not greater than~$d$.
Let $w = w'x_j$ be a~monomial of degree~$d+1$
for $w' = x_1^{i_1}\ldots x_n^{i_n}$, where $i_1+\ldots+i_n = d$.
From the equality
$$
\alpha(w')\alpha_j w
 = R(w')R(x_j)
 = (\alpha(w')+\alpha_j+1)R(w),
$$
we calculate the coefficient~$k$ in $R(w) = kw$ as
\begin{multline*}
k= \frac{\alpha(w')\alpha_j}{\alpha(w')+\alpha_j+1} \\
 = \frac{\alpha_1^{i_1}\ldots\alpha_n^{i_n}\alpha_j}{
 (\alpha_1+1)^{i_1}\ldots(\alpha_n+1)^{i_n}-\alpha_1^{i_1}\ldots\alpha_n^{i_n}}
 \cdot \frac{1}{ \frac{\alpha_1^{i_1}\ldots\alpha_n^{i_n}}
  {(\alpha_1+1)^{i_1}\ldots(\alpha_n+1)^{i_n} -\alpha_1^{i_1}\ldots\alpha_n^{i_n}} +\alpha_j + 1 } \\
 = \frac{\alpha_1^{i_1}\ldots\alpha_n^{i_n}\alpha_j}
 {\alpha_1^{i_1}\ldots\alpha_n^{i_n}
 + (1+\alpha_j)\big( (\alpha_1+1)^{i_1}\ldots(\alpha_n+1)^{i_n} - \alpha_1^{i_1}\ldots\alpha_n^{i_n}\big)} \\
  = \frac{\alpha_1^{i_1}\ldots\alpha_n^{i_n}\alpha_j}{
  (\alpha_1+1)^{i_1}\ldots(\alpha_n+1)^{i_n}(\alpha_j+1)
   - \alpha_1^{i_1}\ldots\alpha_n^{i_n}\alpha_j}
  = \alpha(w),
\end{multline*}
we are done.

{\bf Corollary 4}.\,Let $R$ be an injective monomial RB-operator
of weight\,1 on $\Bbbk[x_1,\ldots,x_n]$. Then $R = -\id$.

{\sc Proof}.
Suppose that $R(1)\not\in \Bbbk$, i.e.,
$R(1) = \alpha x_1^{i_1}\ldots x_n^{i_n}$ with $\alpha\neq0$.
Then, analogous to~\eqref{R(1)R(1)Nonzero}, we get
$R(x_1^{i_1}\ldots x_n^{i_n}) = (1/2)x_1^{i_1}\ldots
x_n^{i_n}(\alpha x_1^{i_1}\ldots x_n^{i_n} - 1)$
is not a~monomial, a~contradiction.
Thus, $R(1) = -1$ and $R$ is splitting.

Since $\ker(R) = (0)$, the ideal $\Bbbk_0[x_1,\ldots,x_n]$ is $R$-invariant.
We may apply Theorem~3. Since $R$ acts on $\Imm(R)$ as $-\id$,
we derive that $\alpha_1 = \ldots = \alpha_n = -1$.
Thus, $R = -\id$. Corollary is proved.

We should say that a direct analogue of Theorem~3 does not hold in the case
of RB-operators of weight zero. By Example~2, we have a collection
$$
S = \{J_a\circ l_f\mid f\in \Bbbk[x]\setminus\{0\},\,a\in \Bbbk\}
$$
of injective RB-operators of weight zero on $\Bbbk[x]$.
If $f$ is a monomial and $a = 0$, then $J_0\circ l_{f}$ is
an injective monomial RB-operator of weight zero on $\Bbbk[x]$.
There is not yet confirmed conjecture about all injective RB-operators on $\Bbbk[x]$.

{\bf Conjecture} (Guo, Rosenkranz, Zheng, 2015~\cite{Monom2}).
The set of all injective RB-operators of weight zero on $\Bbbk[x]$ coincides with~$S$.

{\bf Example 5}.
Let $\alpha_1,\ldots,\alpha_n$ be such nonzero scalars from $\Bbbk$ that
all sums
$\sum\limits_{i_j\in\{1,\ldots,n\}}\frac{1}{\alpha_{i_j}}$ are nonzero.
Define a linear operator $R$ on $\Bbbk_0[x_1,\ldots,x_n]$ as follows:
\begin{equation}\label{MonomialInjectiveOnk[X]Zero}
R(x_1^{i_1}\ldots x_n^{i_n})
 = \frac{x_1^{i_1}\ldots x_n^{i_n}}{\frac{i_1}{\alpha_1}+\ldots+\frac{i_n}{\alpha_n}},
\end{equation}
where $i_1,\ldots,i_n\geq0$, $i_1^2+\ldots+i_n^2>0$.
Then $R$ is an injective monomial RB-operator of weight~0 on $\Bbbk_0[x_1,\ldots,x_n]$.

\section{Gradings by the spectrum of an RB-operator}

It is well-known that for a given finite-dimensional algebra~$A$
we have a grading on~$A$
by the spectrum of its derivation or automorphism~\cite{Jacobson}.
At the end of this work let us state an analogous results for RB-operators,
which are connected to the formulas~\eqref{MonomialInjectiveOnk[X]}
and~\eqref{MonomialInjectiveOnk[X]Zero}.

In this section, a ground field $\Bbbk$ may have a positive characteristic too.

For nonzero scalars $\lambda,\mu\in \Bbbk$, put
$$
\lambda\circ \mu
 = \begin{cases}
 \frac{\lambda\mu}{\lambda+\mu+1}, & \lambda+\mu \neq -1, \\
 \mbox{not defined}, & \mbox{otherwise}.
 \end{cases}
$$
Note that the partially defined operation $\circ$ on $\Bbbk_0$
is commutative and associative in the sense that the equality
$$
(\lambda\circ \mu)\circ \nu
 = \frac{\lambda\mu\nu}{(\lambda+1)(\mu+1)(\nu+1)-\lambda\mu\nu}
 = \lambda\circ(\mu\circ \nu)
$$
holds when all four involved products are defined.

{\bf Example 6}.
The set $\mathbb{R}_{>0}$ of all positive real numbers
under the product $\circ$ is a~semigroup.
Moreover, it is isomorphic to the semigroup
$\langle \mathbb{R}_{>1},\cdot\rangle$.
Indeed, define
$\varphi\colon \mathbb{R}_{>0}\to \mathbb{R}_{>1}$
as $\varphi(x) = 1 + 1/x$. It is a bijection and
$$
\varphi(x\circ y)
 = 1 + \frac{x+y+1}{xy}
 = 1 + \frac{1}{x} + \frac{1}{y} + \frac{1}{xy}
 = \left(1+\frac{1}{x}\right)\left(1+\frac{1}{y}\right)
 = \varphi(x)\varphi(y).
$$

Given a finite-dimensional algebra~$A$ over an algebraically closed field $\Bbbk$
and an RB-operator~$R$ of weight~1 on $A$,
consider the generalized eigenvalue decomposition
\begin{equation}\label{EigenDecomp}
A = \bigoplus\limits_{\lambda\in\Spec(R)} A_{\lambda}.
\end{equation}

{\bf Proposition 2}.
Let $A$ be a finite-dimensional algebra over an algebraically closed field $\Bbbk$
and let $R$ be an RB-operator of weight~1 on $A$.
Let $\lambda,\mu\in\Spec(R)$ be such that $\lambda,\mu\neq0$.
Then
$A_\lambda A_\mu \subseteq
 \begin{cases}
(0), & \lambda+\mu=-1\ or\ \lambda\circ \mu\not\in\Spec(R), \\
A_{\lambda\circ \mu}, & otherwise. \\
\end{cases}$

{\sc Proof}.
Let $u\in A_\lambda$ and $v\in L_\mu$.
So, $(R-\lambda\id)^k u = 0$ and $(R-\mu\id)^l v = 0$.
Let us prove the statement by induction on $k+l$.
Let $k+l = 2$. Then we have
$$
\lambda\mu uv
 = R(u)R(v)
 = R(R(u)v+uR(v)+uv)
 = (\lambda+\mu+1)R(uv).
$$
If $\lambda+\mu = -1$, then $uv = 0$.
Otherwise, $R(uv) = (\lambda\circ\mu) uv$.
From this, we prove the base case.

Suppose that we have proved the statement for all numbers
less than $k+l$. Denote
$\tilde{u} = (R-\lambda\id)u$ and
$\tilde{v} = (R-\mu\id)v$.
So,
\begin{equation}\label{GradIndStep}
\begin{gathered}
R(u)R(v)
 = (\lambda u + \tilde{u})(\mu v + \tilde{v})
 = \lambda\mu uv + \lambda u \tilde{v} + \mu\tilde{u}v + \tilde{u}\tilde{v}, \\
R(R(u)v+uR(v)+uv)
 = (\lambda+\mu+1)R(uv) + R(\tilde{u}v + u\tilde{v}).
\end{gathered}
\end{equation}
If $\lambda+\mu = -1$, then $uv = 0$ since all other products
in~\eqref{GradIndStep} are zero by induction.
Otherwise, we conclude that
$(\lambda+\mu+1)(R-\lambda\circ\mu\id)(uv)\in A_{\lambda\circ\mu}$, i.e.,
$uv\in A_{\lambda\circ\mu}$.
Proposition is proved.

{\bf Example 7}.
Let $A$ be the quotient of $\Bbbk_0[x]$ by the ideal generated by $x^{N+1}$.
Let $\Bbbk$ be either $\mathbb{Q}$ or $\mathbb{Z}_p$ for some prime~$p$.
Then a linear operator $R$ defined on~$A$ as follows
$$
R(x^i) = \dfrac{x^i}{2^i-1},\quad i=1,\ldots,N,
$$
is a monomial RB-operator of weight~1 on~$A$.
To avoid division by zero in the case $\Bbbk = \mathbb{Z}_p$,
we require that none of the numbers
$2^2-1$, $2^3-1$, \ldots, $2^N-1$ is divided by~$p$.
For this, we may restrict $N<\log_2(p+1)$.
Multiplying the unit $i$ times, we get
$1^i_\circ = \dfrac{1}{2^i-1}$.
So, we have the decomposition
$$
A = A_1\oplus A_{1^2_\circ}\oplus\ldots\oplus A_{1^N_\circ},\quad A_{1^i_\circ} = \Span\{x^i\}.
$$
Let us check that Proposition~2 holds true in this example.
When $i+j\leq N$, we have
$x^i\cdot x^j = x^{i+j}\in A_{1^{i+j}_\circ}$, since
$1^i_\circ \circ 1^j_\circ = 1^{i+j}_\circ$.
When $i+j>N$, we have $x^i\cdot x^j = 0$ and
$1^{i+j}_\circ\not\in\Spec(R)$.

Let us give a particular example of this construction
for $N = 3$ and $p = 5$.
In this case $2^2-1 = 3$ and $2^3-1 = 7$ are not divided by~5.
So, we have $R(x) = x$, $R(x^2) = 2x^2$, $R(x^3) = 3x^3$ and a decomposition
$A = A_1\oplus A_2\oplus A_3$ with $A_i = \Span\{x^i\}$.
By Proposition~2,
$x\cdot x^2 = x^3\in A_3$, since
$1\circ 2 = 3$.
On the other hand, $1\circ 3$ is not defined,
and we have $x\cdot x^3 = 0$ correspondingly.
Finally, $2\circ 3 = 1$, and $x^2\cdot x^3 = 0\in A_1$.

Now, let us proceed to the case of weight zero.
For nonzero scalars $\lambda,\mu\in \Bbbk$, put
$$
\lambda *\mu
 = \begin{cases}
 \frac{\lambda\mu}{\lambda+\mu}, & \lambda\neq -\mu, \\
 \mbox{not defined}, & \mbox{otherwise}.
 \end{cases}
$$
The operation $*$ on $\Bbbk_0$ is again commutative and associative, i.e.,
$$
(\lambda*\mu)*\nu
 = \frac{\lambda\mu\nu}{\lambda\mu+\lambda\nu+\mu\nu}
 = \lambda*(\mu *\nu)
$$
holds when all products involved are defined.

{\bf Example 8}.
The set $\mathbb{R}_{>0}$ of all positive real numbers
under the product $\circ$ is a~semigroup, and
it is isomorphic to the semigroup
$\langle \mathbb{R}_{>0},+\rangle$.
Indeed, define
$\varphi\colon \mathbb{R}_{>0}\to \mathbb{R}_{>0}$
as $\varphi(x) = 1/x$. It is a bijection and
$\varphi(x*y)
 = \frac{x+y}{xy}
 = \frac{1}{x} + \frac{1}{y}
 = \varphi(x) + \varphi(y)$.

We again consider the generalized eigenvalue decomposition~\eqref{EigenDecomp}
of an algebra~$A$ by an RB-operator~$R$ of weight zero.

{\bf Proposition 3}.
Let $A$ be a finite-dimensional algebra over
an algebraically closed field $\Bbbk$
and let $R$ be an RB-operator of weight~0 on $A$.
Let $\lambda,\mu\in\Spec(R)$ be such that $\lambda,\mu\neq0$.
Then
$A_\lambda A_\mu \subseteq
 \begin{cases}
(0), & \lambda+\mu=0\ or\ \lambda*\mu\not\in\Spec(R), \\
A_{\lambda*\mu}, & otherwise. \\
\end{cases}$

{\sc Proof}.
Analogous to the proof of Proposition~2.

{\bf Example 9}.
We deal with the same algebra as in Example~7.
A linear operator $R$ defined on~$A$ as follows
$$
R(x^i) = x^i/i,\quad i=1,\ldots,N,
$$
is a monomial RB-operator of weight~0 on~$A$. Here we need $N<p$.
Multiplying the unit $i$ times, we get
$1^i_* = 1/i$. So, we have the decomposition
$$
A = A_1\oplus A_{1/2}\oplus\ldots\oplus A_{1/N},\quad A_{1/i} = \Span\{x^i\}.
$$
When $i+j\leq N$, we have
$x^i\cdot x^j = x^{i+j}\in A_{1/(i+j)}$, since
$1^i_{*} * 1^j_{*}
 = \frac{1}{i}*\frac{1}{j}
 = \frac{1}{i+j}
 = 1^{i+j}_*$.
When $i+j>N$, we have $x^i\cdot x^j = 0$ and
$1^{i+j}_*\not\in\Spec(R)$.

Let us consider Example~9 with $N = 3$ and $p = 5$.
Then $R(x) = x$, $R(x^2) = 3x^2$, $R(x^3) = 2x^3$,
and $A = A_1\oplus A_2\oplus A_3$.
Here $x^2\cdot x^3 = 0$, since $(1/2)*(1/3)$ is not defined.
Further, $x\cdot x^2 = x^3\in A_2$, since $1*3 = 2$.
Finally, $x\cdot x^3 = 0$ as $1*2 = 4\not\in\Spec(R)$.

{\bf Remark 3}.
Note that Example~7 is nothing more as a quotient of the RB-algebra~\eqref{RBNonzeroNonunitary}
with $\alpha = 1$ by the Rota---Baxter-ideal generated by $x^{N+1}$.
Analogously, the quotient of the RB-algebra from Example~3 gives
the RB-operator from Example~9.

{\bf Remark 4}.
It is important that we exclude unit in both Examples~7 and~9, it means
we do not consider the quotient $B$ of the polynomial algebra $\Bbbk[x]$
by the ideal generated by $x^{N+1}$.
Indeed, in this case spectrum of every RB-operator of weight~$\lambda$
on~$B$ is a subset of $\{0,-\lambda\}$~\cite{Unital,Spectrum}.

The results of Propositions~2 and~3 can be useful in attempt to construct
universal enveloping associative Rota---Baxter algebra for a given Lie Rota---Baxter algebra,
see the exact formulation of the problem in~\cite{postLiepostAs}.

\section*{Acknowledgements}

Author was supported by the Program of fundamental scientific researches
of the Siberian Branch of Russian Academy of Sciences, I.1.1, project 0314-2019-0001.
Author is grateful to Korean Institute of Advanced Study for support and hospitality
(August 2019). Author thanks an anonymous reviewer for the helpful remarks.

\medskip

\noindent Vsevolod Gubarev \\
Sobolev Institute of Mathematics \\
Acad. Koptyug ave. 4, 630090 Novosibirsk, Russia \\
Novosibirsk State University \\
Pirogova str. 2, 630090 Novosibirsk, Russia \\
e-mail: wsewolod89@gmail.com

\begin{thebibliography}{67}
\bibitem{Aguiar}
M. Aguiar.
Pre-Poisson algebras, Lett. Math. Phys. {\bf 54} (2000) 263--277.

\bibitem{FardThesis}
K. Ebrahimi-Fard.
Rota-Baxter Algebras and the Hopf Algebra of Renormalization,
Ph.D. Thesis, University of Bonn, 2006.

\bibitem{Baxter}
G. Baxter.
An analytic problem whose solution follows from a simple algebraic identity. Pacific J. Math. {\bf 10} (1960) 731--742.

\bibitem{BGP}
P. Benito, V. Gubarev, A. Pozhidaev.
Rota---Baxter operators on quadratic algebras, Mediterr. J. Math. {\bf 15} (2018), 23 p. (N189).

\bibitem{Unital}
V. Gubarev.
Rota---Baxter operators on unital algebras, Mosc. Math. J. (accepted), arXiv.1805.00723v3, 43~p.

\bibitem{Spectrum}
V. Gubarev.
Spectrum of Rota---Baxter operators, preprint, 20~p.

\bibitem{postLiepostAs}
V. Gubarev.
Embedding of post-Lie algebras into postassociative algebras, 
New Trends in Algebra and Combinatorics. Proc. of the 3rd Intern. 
Congress in Algebra and Combinatorics 
(Ed. by K.P. Shum et al), World Scientific, p. 57--67, 2020.

\bibitem{GuoMonograph}
L. Guo.
An Introduction to Rota---Baxter Algebra. Surveys of Modern Mathematics, vol. 4,
Intern. Press, Somerville (MA, USA); Higher education press, Beijing, 2012.

\bibitem{GuoLiu}
L. Guo, Z. Liu.
Rota-Baxter operators on generalized power series rings, J. Algebra Appl. {\bf 8} (2009) 557--564.

\bibitem{Monom2}
L. Guo, M. Rosenkranz, and S.H. Zheng.
Rota-Baxter operators on the polynomial algebras,
integration and averaging operators,
Pacific J. Math. (2) {\bf 275} (2015) 481--507.

\bibitem{Jacobson}
N. Jacobson.
Lie  algebras. Interscience Tracts on Pure and Applied Mathematics,  no. 10.
Interscience Publishers, N.-Y., 1962.

\bibitem{Mazurek}
R. Mazurek.
Rota---Baxter operators on skew generalized power series rings,
J. Algebra Appl. (7) {\bf 13} (2014), 1650087, 10~p.

\bibitem{Monom}
H. Yu.
Classification of monomial Rota-Baxter operators on $k[x]$,
J. Algebra Appl. {\bf 15} (2016), 1650087, 16~p.

\bibitem{Zhelyabin}
V.N. Zhelyabin.
Jordan bialgebras of symmetric elements and Lie bialgebras,
Siberian Mat. J. (2) {\bf 39} (1998) 261--276.
\end{thebibliography}
\end{document}